\documentclass[12pt,leqno]{amsart}
\newtheorem{thm}{Theorem}[section]
\newtheorem{defi}{Definition}[section]

\newtheorem{rem}{Remark}[section]
\newtheorem{exm}{Example}[section]
\newcommand{\be}{\begin{equation}}
\newcommand{\ee}{\end{equation}}
\newcommand{\bea}{\begin{eqnarray}}
\newcommand{\eea}{\end{eqnarray}}
\newcommand{\beb}{\begin{eqnarray*}}
\newcommand{\eeb}{\end{eqnarray*}}
\usepackage{amssymb,amsfonts,amsthm,setspace,indentfirst,mathrsfs}
\usepackage{minibox}
\usepackage{enumitem}
\numberwithin{equation}{section}
\begin{document}
%
\title[Curvature properties of Robinson-Trautman metric]{\bf{Curvature properties of Robinson-Trautman metric}}
\author[Absos Ali Shaikh, Musavvir Ali and Zafar Ahsan]{Absos Ali Shaikh$^{*1}$, Musavvir Ali$^2$ Zafar Ahsan$^3$ }
\date{\today}
\address{\noindent$^{1 \& 2}$ Department of Mathematics,
\newline Department of Mathematics
\newline Aligarh Muslim University,
\newline Aligarh-202002,
\newline Uttar Pradesh, India}
\email{aask2003@yahoo.co.in, aashaikh@math.buruniv.ac.in}
\email{musavvirali.maths@amu.ac.in}
\address{\noindent$^3$ Faculty of Science and Technology,
\newline University of Islamic Sciences,
\newline Nilai,Malaysia}
\email{zafar.ahsan@radiffmail.com}
\dedicatory{}
\begin{abstract}
The  curvature properties of Robinson-Trautman metric have been investigated. It is shown that Robinson-Trautman metric admits several kinds of pseudosymmetric type structures such as Weyl pseudosymmetric, Ricci pseudosymmetric, pseudosymmetric Weyl conformal curvature tensor etc. Also it is shown that the difference $R\cdot R - Q(S,R)$ is linearly dependent with $Q(g,C)$ but the metric is not Ricci generalized pseudosymmetric. Moreover, it is proved that this metric is Roter type, 2-quasi-Einstein, Ricci tensor is Riemann compatible and its Weyl conformal curvature 2-forms are recurrent. It is also shown that the energy momentum tensor of the metric is pseudosymmetric and the conditions under which such tensor is of Codazzi type and cyclic parallel have been investigated. Finally, we have made a comparison between the curvature properties of Robinson-Trautman metric and Som-Raychaudhuri metric.
\end{abstract}
%
\subjclass[2010]{53B20, 53B25, 53B30, 53B50, 53C15, 53C25, 53C35, 83C15}
\keywords{Robinson-Trautman metric, Einstein field equation, Weyl conformal curvature tensor, pseudosymmetric type curvature condition, 2-quasi-Einstein manifold, Som-Raychaudhuri spacetime}
\maketitle
%

\section{\bf Introduction}\label{intro}
Let $M$ be a connected smooth semi-Riemannian manifold of dimension $n$ $(\geq 3)$ endowed with the semi-Riemannian metric $g$  with signature $(p, n-p)$ and the Levi-Civita connection $\nabla$. If $p = 1$ or $n-1$, then $M$ is Lorentzian and if $p = 0$ or $n$, then $M$ is Riemannian. Let  $R$, $S$ and $\kappa $ be respectively the Riemann-Christoffel curvature tensor of type (0,4), the Ricci tensor of type (0,2) and the scalar curvature of $M$. As a generalization of manifold of constant curvature, Cartan \cite{Cart26} introduced the notion of locally symmetric manifolds defined as $\nabla R =0$. Again, Cartan \cite{Cart46} generalized the concept of locally symmetric manifolds and introduced the notion of semisymmetric manifolds which were latter classified by Szab\'o (\cite{Szab82}, \cite{Szab84}, \cite{Szab85}). During the study of totally umbilical submanifolds of semisymmetric manifolds, Adam\'{o}w and Deszcz \cite{AD83} introduced the notion of pseudosymmetric manifolds which was latter on studied by many authors (cf. \cite{Desz92}, \cite{Desz93}, \cite{DGHS11}, \cite{DG98}, \cite{DHV08}, \cite{HV04}, \cite{HV07a}, \cite{HV07b}, \cite{SDHJK15} and also references therein). We note that the recent trend of modern mathematical research is the abstraction, generalization or extension, existence, characterization, classification and finally applications. We mention that during the last four decades generalized concept of pseudosymmetry by Deszcz has been characterized and classified with existence in several papers by Deszcz and his coauthors (see \cite{DD91a}, \cite{DD91b}, \cite{Desz93}, \cite{DGHZ15}, \cite{SDHJK15} and also references therein); and the pseudosymmetric spacetimes were studied by Deszcz and his co-authors (see, \cite{DDVV94}, \cite{DHV04}, \cite{DVV91}). It may be noted that G\"{o}del spacetime (\cite{DHJKS14}, \cite{Gode49}), Som-Raychaudhuri spacetime (\cite{SK16srs}, \cite{SR68}), Reissner-Nordstr\"{o}m spacetime \cite{Kowa06} and Robertson-Walker spacetimes (\cite{ADEHM14}, \cite{DDHKS00}, \cite{DK99}), are models of different pseudosymmetric type structures.\\
\indent It is known that the Robinson-Trautman metrics are the spacetimes which admit a geodesic, hypersurface orthogonal, shear-free and expanding null congruences. These metrics are the general relativistic analogues of Lienard-Wiechert solutions of Maxwell equations. Moreover, the Robinson-Trautman spacetimes generalize the assumption of spherical symmetry by having topologically equivalent two-spheres rather than strictly two spheres. From the physical point of view, the Robinson-Trautman metrics can be thought of as representing an isolated gravitationally radiating system. The Robinson-Trautman metric of Petrov type II is given by (\cite{RT62}, \cite{ste03}, \cite{zun75})
$$
ds^2 = -2(U^0 - 2 \gamma^0 r - \Psi^0_2 r^{-1})du^2 + 2 du dr - \frac{r^2}{2 P^2} d\zeta d\bar\zeta,
$$
where $U^0, \gamma^0, \Psi^0_2$ are constants and $P$ is a nowhere vanishing function of $\zeta$ and $\bar\zeta$. For the sake of simplicity, we write these constants $U^0, \gamma^0, \Psi^0_2$ as $a, b, q$ respectively, $P$ as $f$ and the coordinate $u$ as $t$. The above metric can be written as
\be\label{rtm}
ds^2 = -2(a - 2 b r - q r^{-1})dt^2 + 2 dt dr - \frac{2 r^2}{f^2} d\zeta d\bar\zeta.
\ee
It is easy to check that this metric can be written as a warped product metric (\cite{BO69}, \cite{Kruc57}) given by
$$ds^2 = d\bar s^2 + r^2 d\tilde s^2,$$
where $d\bar s^2 = -2(a - 2 b r - q r^{-1})dt^2 + 2 dt dr$ and $d\tilde s^2 = - \frac{2}{f^2} d\zeta d\bar\zeta$.\\
\indent The main object of the present paper is to investigate the pseudosymmetric type structures admitting by the Robinson-Trautman metric \eqref{rtm}. Lanczos potentials of the metric are studied by Ahsan and Bilal \cite{AZBM10}. In the study of differential geometry, pseudosymmetric structures play an important role due their applications in relativity and cosmology (see, \cite{DDVV94}, \cite{DHV04}, \cite{DHJKS14}, \cite{DVV91}, \cite{SK16srs}, \cite{SM17} and also references therein). It is noteworthy to mention that the present paper is unconventional as by considering a physically relevant metric, viz., Robinson-Trautman metric, we would like to investigate the curvature restricted geometric structures admitted by such metric. Then Robison-Trautman spacetime can be considered as a model of such pseudosymmetric structures.\\
\indent The paper is organized as follows. Section 2 deals with the different pseudosymmetric structures which are essential to investigate the curvature properties of the Robinson-Trautman metric. Section 3 is concerned with different components of several curvature tensors which are useful to investigate the pseudosymmetric type structures admitting by this metric. Section 4 is devoted to the conclusion as theorems and remarks. Section 5 deals with the investigation under which the energy momentum tensor is covariantly constant, Codazzi type and cyclic parallel. In the last section we have made a comparison (similarities and dissimilarities) between the curvature properties of Robinson-Trautman metric and Som-Raychaudhuri metric. It is shown that both the metric are 2-quasi-Einstein and pseudosymmetric Weyl conformal curvature tensor, whereas Robinson-Trautman metric is Deszcz pseudosymmetric but Som-Raychaudhuri metric is Ricci generalized pseudosymmetric.
\section{\bf Curvature restricted geometric structures}
 It is well known that a curvature restricted geometric structure is a geometric structure on $M$ obtained by imposing a restriction on its curvature tensors by means of covariant derivatives of first order or higher orders. We will now explain some useful notations and definitions of various curvature restricted geometric structures. For two symmetric $(0,2)$-tensors $A$ and $E$, their Kulkarni-Nomizu product $A\wedge E$ is defined as (see e.g. \cite{DGHS11}, \cite{Glog02}):
\begin{eqnarray*}
(A\wedge E)(X_1,X_2,X_3,X_4) &=& A(X_1,X_4)E(X_2,X_3) + A(X_2,X_3)E(X_1,X_4)\\
&-& A(X_1,X_3)E(X_2,X_4) - A(X_2,X_4)E(X_1,X_3),
\end{eqnarray*}
where $X_1,X_2,X_3,X_4 \in \chi(M)$, the Lie algebra of all smooth vector fields on $M$. Throughout the paper we will consider $X, Y, X_1, X_2, \cdots \in \chi(M)$.\\
\indent Again as an invariant under certain class of transformations on a semi-Riemannian manifold $M$, there arise various curvature tensors on $M$, such as conformal curvature tensor $C$, projective curvature tensor $P$, concircular curvature tensor $W$ and conharmonic curvature tensor $K$ (\cite{Ishi57}, \cite{YK89}) etc. In terms of Kulkarni-Nomizu product $C$, $W$, $K$ and the Gaussian curvature tensor $\mathfrak G$ are given by
$$C = R-\frac{1}{n-2}(g\wedge S) + \frac{r}{2(n-2)(n-1)}(g\wedge g),$$
$$W = R-\frac{r}{2n(n-1)}(g\wedge g),$$
$$K = R-\frac{1}{n-2}(g\wedge S) \ \mbox{ and } \ \mathfrak G = \frac{1}{2}(g\wedge g).$$
The projective curvature tensor $P$ is given by
$$
P(X_1, X_2, X_3, X_4) = R(X_1, X_2, X_3, X_4) - \frac{1}{n-1}[g(X_1, X_4)S(X_2, X_3)-g(X_2, X_4)S(X_1, X_3)].
$$
A $(0,4)$ tensor is called a generalized curvature tensor if it posses the symmetry like $R$, i.e., a $(0,4)$ tensor $D$ is said to be a generalized curvature tensor (\cite{DGHS11}, \cite{SDHJK15}, \cite{SK14}) if
$$D(X_1,X_2,X_3,X_4)+D(X_2,X_3,X_1,X_4)+D(X_3,X_1,X_2,X_4)=0,$$
$$D(X_1,X_2,X_3,X_4)+D(X_2,X_1,X_3,X_4)=0 \ \ \mbox{and}$$
$$D(X_1,X_2,X_3,X_4)=D(X_3,X_4,X_1,X_2).$$
We note that if $A$ and $E$ are both symmetric $(0,2)$-tensors, then $A\wedge E$ is a generalized curvature tensor, and linear combination of two generalized curvature tensors is also a generalized curvature tensor. Hence $C$, $W$ and $K$ are generalized curvature tensors but $P$ is not a generalized curvature tensor.\\
\indent Now for a $(0,4)$-tensor $D$ and a symmetric $(0,2)$-tensor $E$ we can define three endomorphisms $\mathcal E$, $\mathscr{D}(X,Y)$ and $X\wedge_E Y$ as (\cite{DDHKS00}, \cite{DGHS11})
$$
E(X, Y) = g(\mathcal E X, Y), \ \mathscr{D}(X,Y)X_1 = \mathcal D(X,Y)X_1 \ \mbox{and} \ (X\wedge_A Y)X_1 = A(Y,X_1)X-A(X,X_1)Y
$$
respectively, where $\mathcal D$ is the corresponding $(1,3)$-tensor of $D$, i.e.,
$$D(X_1,X_2,X_3,X_4) = g(\mathcal D(X_1,X_2)X_3, X_4).$$
\indent One can now easily operate $\mathscr{D}(X,Y)$ and $X\wedge_A Y$ on a $(0,k)$-tensor $H$, $k\geq 1$, and obtain two $(0,k+2)$-tensors $D\cdot H$ and $Q(A,H)$ respectively as follows (see \cite{DG02}, \cite{DGHS98}, \cite{DH03}, \cite{SDHJK15}, \cite{SK14}, \cite{SKgrt}, \cite{Tach74} and also references therein):
$$D\cdot H(X_1,X_2,\cdots,X_k,X,Y) = -H(\mathcal D(X,Y)X_1,X_2,\cdots,X_k) - \cdots - H(X_1,X_2,\cdots,\mathcal D(X,Y)X_k).$$
and
\beb
&&Q(A,H)(X_1,X_2, \ldots ,X_k,X,Y) = ((X \wedge_A Y)\cdot H)(X_1,X_2, \ldots ,X_k)\\
&&= A(X, X_1) H(Y,X_2,\cdots,X_k) + \cdots + A(X, X_k) H(X_1,X_2,\cdots,Y)\\
&& - A(Y, X_1) H(X,X_2,\cdots,X_k) - \cdots - A(Y, X_k) H(X_1,X_2,\cdots,X).
\eeb
In terms of local coordinates system, $D\cdot H$ and $Q(A,H)$ \cite{Desz03a} can be written as
\beb
(D\cdot H)_{i_1 i_2\cdots i_k j l} &=& -g^{pq}\left[D_{jli_1 q}H_{p i_2\cdots i_k} + \cdots + D_{jli_k q}H_{i_1 i_2\cdots p}\right],
\eeb
\beb
Q(A, H)_{i_1 i_2\cdots i_k j l} &=& A_{li_1}H_{j i_2\cdots i_k} + \cdots + A_{li_k}H_{i_1 i_2\cdots j}\\
																	& & -A_{ji_1}H_{l i_2\cdots i_k} - \cdots - A_{ji_k}H_{i_1 i_2\cdots l}.
\eeb
\begin{defi}$($\cite{AD83}, \cite{Cart46}, \cite{Desz92}, \cite{SK14}, \cite{SKppsn}, \cite{SKppsnw}, \cite{SRK15}, \cite{Szab82}$)$
A semi-Riemannian manifold $M$ is said to be $H$-semisymmetric type if $D\cdot H = 0$ and it is said to be $H$-pseudosymmetric type if $\left(\sum\limits_{i=1}^k c_i D_i\right)\cdot H = 0$ for some scalars $c_i$'s, where $D$ and each $D_i$, $i=1,\ldots, k$, $(k\ge 2)$, are (0,4) curvature tensors. In particular, if $c_i$'s are all constants, then the manifold is called $H$-pseudosymmetric type manifold of constant type or otherwise non-constant type.
\end{defi}
\indent In particular, if $D = R$ and $T=R$ (resp., $S$, $C$, $W$ and $K$), then $M$ is called semisymmetric (resp., Ricci, conformally, concircularly and conharmonically semisymmetric). Again, if $i =2$, $D_1 = R$, $D_2 = \mathfrak G$ and $D= R$ (resp., $S$, $C$, $W$ and $K$), then $M$ is called Deszcz pseudosymmetric (resp., Ricci, conformally, concircularly and conharmonically pseudosymmetric). Especially, if $i =2$, $D_1 = C$, $D_2 = \mathfrak G$ and $D =C$, then $M$ is called a manifold of pseudosymmetric Weyl conformal curvature tensor. Again, if $i =2$, $D_1 = R$, $D_2 = \wedge_S$ and $D = R$, then $M$ is called Ricci generalized pseudosymmetric \cite{DD91a}.\\
\indent A semi-Riemannian manifold $M$ is said to be Einstein if its Ricci tensor is a scalar multiple of the metric tensor $g$. We note that in this case $S=\frac{\kappa}{n}g$ and $\kappa$ is scalar curvature. Again $M$ is called quasi-Einstein if at each point of $M$, rank of $(S - \alpha g)$ is less or equal to 1 for a scalar $\alpha$. Thus in this case $S = \alpha g + \beta \Pi \otimes \Pi$, where $\alpha$ and $\beta$ are scalars and $\Pi$ is an 1-form. In particular, if $\alpha = 0$, then a quasi-Einstein manifold is called Ricci simple. We note that G\"{o}del spacetime is Ricci simple \cite{DHJKS14} and perfect fluid spacetimes are quasi-Einstein (see, \cite{SB05}, \cite{SB06}). For quasi-Einstein spacetimes we refer the reader to see \cite{SYH09} and also references therein.
\begin{defi}
A semi-Riemannian manifold is said to be 2-quasi-Einstein if at each point of $M$, rank of $(S - \alpha g)$ is less or equal to 2 for a scalar $\alpha$.
\end{defi}
\indent It may be noted that Som-Raychaudhuri spacetime is 2-quasi-Einstein but not quasi-Einstein \cite{SK16srs}. For recent results on 2-quasi-Einstein manifolds, we refer the reader to see \cite{DGHZ-2016}, \cite{DGJZ-2016} and \cite{DGP-TV-2015} .
\begin{defi} \cite{Shai09}
A semi-Riemannian manifold $M$ is said to be pseudo quasi-Einstein if
\be\label{pqe}
S = \alpha g + \beta \Pi \otimes \Pi + \gamma E
\ee
holds on $U$ for some scalars $\alpha, \beta, \gamma$, an 1-form $\Pi$ and a trace free $(0,2)$-tensor $E$ such that the associated vector field $V$ corresponding to $\Pi$ satisfies $E(X,V) =0$.
\end{defi}
\indent In particular, if $E = \Phi \otimes \Phi$ (resp., $\Pi \otimes \Phi + \Phi \otimes \Pi$) then a pseudo quasi-Einstein manifold reduces to a generalized quasi-Einstein manifold by De and Ghosh \cite{DG04} (resp., by Chaki \cite{Chak01}).\\
\indent To get the decomposition of Riemann-Christoffel curvature tensor $R$, Deszcz (\cite{Desz03}, \cite{Desz03a}) introduced the notion of Roter type manifold, which was latter extended as generalized Roter type manifold \cite{SDHJK15}.
\begin{defi}
A semi-Riemannian manifold $M$ is said to be a Roter type manifold $($\cite{Desz03}, \cite{Desz03a}, \cite{DGHS11}, \cite{DGP-TV-2011}, \cite{DPSch-2013} and \cite{Glog-2007}$)$ $($resp., generalized Roter type manifold $($\cite{DGJPZ13}, \cite{DGJZ-2016}, \cite{DGP-TV-2015}, \cite{SDHJK15}, \cite{SKgrt}, \cite{SK16} and \cite{Saw-2015}$))$ if
$$R = N_1 (g\wedge g) + N_2 (g\wedge S) + N_3 (S\wedge S) \ \mbox{ and}$$
$$\left[\mbox{resp., } R = L_1 (g\wedge g) + L_2 (g\wedge S) + L_3 (S\wedge S) + L_4 (g\wedge S^2) + L_5 (S\wedge S^2) + L_6 (S^2\wedge S^2)\right]$$
holds for some $N_i, L_j \in C^{\infty}(M)$, $1\le i\le 3$ and $1\le j\le 6$.
\end{defi}
\begin{defi}(\cite{Bess87}, \cite{SKgrt})
A semi-Riemannian manifold $M$ is said to be $Ein(2)$ if $S^2$, $S$ and $g$ are linearly dependent, where $S^2$ is the second level Ricci tensor defined as $S^2(X,Y) = S(\mathcal S X, Y)$.
\end{defi}
\indent Between the class of Ricci symmetric manifolds and the class of manifolds of constant scalar curvature, Gray \cite{Gray78} presented two new classes of manifolds, namely, manifolds of Codazzi type Ricci tensor and manifolds of cyclic parallel Ricci tensor.
\begin{defi}
A semi-Riemannian manifold $M$ is said to be of Codazzi type Ricci tensor (resp. cyclic parallel Ricci tensor) (see, \cite{DHJKS14}, \cite{Gray78} and references therein) if
\[(\nabla_{X_1} S)(X_2, X_3) = (\nabla_{X_2} S)(X_1, X_3)\]
\[\left(\mbox{resp.} \ \ \ (\nabla_{X_1} S)(X_2, X_3) + (\nabla_{X_2} S)(X_3, X_1) + (\nabla_{X_3} S)(X_1, X_2) = 0 \ \right).\]
\end{defi}
Recently Mantica and Molinari showed \cite{MM12a} that the sufficient condition of Codazzi type for Derdzinski and Shen theorem \cite{DS83} can be replaced by a more general condition, and they introduced the notion of Riemann-compatible tensor along with  compatible tensors for other curvature tensors.
\begin{defi}\label{def2.8}
Let $D$ be a $(0,4)$-tensor and $E$ be a symmetric $(0, 2)$-tensor on $M$. Then $E$ is said to be $D$-compatible (\cite{DGJPZ13}, \cite{MM12b}, \cite{MM13}) if
\[
D(\mathcal E X_1, X,X_2,X_3) + D(\mathcal E X_2, X,X_3,X_1) + D(\mathcal E X_3, X,X_1,X_2) = 0
\]
holds. Again an 1-form $\Pi$ is said to be $D$-compatible if $\Pi\otimes \Pi$ is $D$-compatible.
\end{defi}
\indent Generalizing the concept of recurrent manifold (\cite{Ruse46}, \cite{Ruse49a}, \cite{Ruse49b}, \cite{Walk50}), recently Shaikh et al. \cite{SRK16} introduced the notion of super generalized recurrent manifold along with its characterization and existence by proper example.
\begin{defi}
A semi-Riemannian manifold $M$ is said to be super generalized recurrent manifold (\cite{SK14}, \cite{SRK16}, \cite{SKA16}) if
$$
\nabla R = \Pi \otimes R + \Phi \otimes (S\wedge S) + \Psi \otimes (g\wedge S) + \Theta \otimes (g\wedge g)
$$
holds on $\{x\in M: R \neq 0 \mbox{ and any one of } S\wedge S, g\wedge S \mbox{ is non-zero at $x$}\}$ for some 1-forms $\Pi$, $\Phi$, $\Psi$ and $\Theta$, called the associated 1-forms. Especially, if $\Phi = \Psi = \Theta = 0$ (resp., $\Psi = \Theta = 0$ and $\Phi = \Theta = 0$), then the manifold is called recurrent (\cite{Ruse46}, \cite{Ruse49a}, \cite{Ruse49b}, \cite{Walk50}) (resp., weakly generalized recurrent (\cite{SR11}, \cite{SAR13}) and hyper generalized recurrent (\cite{SP10}, \cite{SRK15})) manifold.
\end{defi}
\indent Again as a generalization of locally symmetric manifold and recurrent manifold, Tam$\acute{\mbox{a}}$ssy and Binh \cite{TB89} introduced the notion of weakly symmetric manifolds.
\begin{defi}
Let $D$ be a (0, 4)-tensor on a semi-Riemannian manifold $M$. Then $M$ is said to be weakly $D$-symmetric (\cite{TB89}, \cite{SK12}) if
\beb
&&\nabla_X  D(X_1, X_2, X_3, X_4) = \Pi(X) D(X_1, X_2, X_3, X_4) + \Phi(X_1) D(X_1, X_2, X_3, X_4)\\
&& \hspace{2cm} + \overline \Phi(X_2) D(X_1, X, X_3, X_4) + \Psi(X_3) D(X_1, X_2, X, X_4) + \overline \Psi(X_4) D(X_1, X_2, X_3, X)
\eeb
holds $\forall~ X, X_i \in \chi(M)$ $(i =1,2,3,4)$ and some 1-forms $\Pi, \Phi, \overline \Phi, \Psi$ and $\overline \Psi$  on $\{x\in M: R_x \neq 0\}$. In particular, if $\frac{1}{2}\Pi = \Phi = \overline \Phi = \Psi = \overline \Psi$, then the manifold is called Chaki $D$-pseudosymmetric manifold \cite{Chak87}.
\end{defi}
\indent It is noteworthy to mention that if $D$ is a generalized curvature tensor then the above defining condition of weak symmetry holds for $\frac{\Phi + \Psi}{2}$ in lieu of $\Phi, \overline \Phi, \Psi$ and $\overline \Psi$. It is also noted that the notion of Chaki pseudosymmetry is different from Deszcz pseudosymmetry. For details about the defining condition of weak symmetry and the interrelation between weak symmetry and Deszcz pseudosymmetry, we refer the reader to see \cite{SDHJK15} and also references therein.\\
\indent For a generalized curvature tensor $D$ and a symmetric (0, 2)-curvature tensor $Z$, one can define the curvature 2-form $\Omega_{(D)l}^m$ (\cite{Bess87}, \cite{LR89}) for $D$ and an 1-form $\Lambda_{(Z)l}$ \cite{SKP03} as
$$\Omega_{(D)l}^m= D_{jkl}^m dx^j \wedge dx^k \mbox{ and } \Lambda_{(Z)l} = Z_{lm} dx^m,$$
where $\wedge$ indicates the exterior product. Now $\Omega_{(D)l}^m$ (resp., $\Lambda_{(Z)l}$) are recurrent if
$$\mathcal D \Omega_{(R)l}^m  = \Pi \wedge \Omega_{(R)l}^m \mbox{ (resp., } \mathcal D \Lambda_{(S)l}  = \Pi \wedge \Lambda_{(S)l}),$$
where $\mathcal D$ is the exterior derivative and $\Pi$ is the associated $1$-form. Recently Mantica and Suh (\cite{MS12a}, \cite{MS13a}, \cite{MS14}) showed that $\Omega_{(D)l}^m$ are recurrent if and only if
\beb\label{man}
&&(\nabla_{X_1} D)(X_2,X_3,X,Y)+(\nabla_{X_2} D)(X_3,X_1,X,Y)+(\nabla_{X_3} D)(X_1,X_2,X,Y) =\\
&&\hspace{1in} \Pi(X_1) D(X_2,X_3,X,Y) + \Pi(X_2) D(X_3,X_1,X,Y)+ \Pi(X_3) D(X_1,X_2,X,Y)
\eeb
and $\Lambda_{(Z)l}$ are recurrent if and only if
$$(\nabla_{X_1} Z)(X_2,X) - (\nabla_{X_2} Z)(X_1,X) = \Pi(X_1) Z(X_2,X) - \Pi(X_2) Z(X_1,X)$$
for an 1-form $\Pi$.
\section{\bf Curvature properties of Robinson-Trautman metric}\label{com}
Setting $\zeta = x^3 + i x^4$, the metric \eqref{rtm} can be written as
\be\label{RTII}
ds^2 = -2(a - 2 b r - q r^{-1})dt^2 + 2 dt dr - \frac{r^2}{f^2} [(dx^3)^2+(dx^4)^2],
\ee
where $f$ is a function of the real variables $x^3$ and $x^4$.
Hence the metric tensor in coordinates $(t,r,x^3,x^4)$ is given by 
$$g = \left(
\begin{array}{cccc}
 -2(a - 2 b r - q r^{-1}) & 1 & 0 & 0 \\
 1 & 0 & 0 & 0 \\
 0 & 0 & - \frac{r^2}{f^2} & 0 \\
 0 & 0 & 0 & - \frac{r^2}{f^2}
\end{array}
\right).$$ 
Again $g$ can be expressed as warped product $\bar g \times_r \tilde g$ with warping function $r$, where base metric $\bar g$ and fiber metric $\tilde g$ are given by
$$\bar g = \left(
\begin{array}{cc}
 -2(a - 2 b r - q r^{-1}) & 1 \\
 1 & 0 
\end{array}
\right), \ \ 
\tilde g = \left(
\begin{array}{cc}
 - \frac{1}{f^2} & 0 \\
 0 & - \frac{1}{f^2}
\end{array}
\right).$$
From above it is clear that the metric is Lorentzian with signature $\{-,+,+,+\}$.\\
\indent The non-zero components of its Riemann-Christoffel curvature tensor $R$, Ricci tensor $S$ and scalar curvature $\kappa$ are given by
$$R_{1212}=-\frac{2 q}{r^3}, \ R_{1313}= R_{1414}=-\frac{2 \left(2 b r^2-q\right) \left(-a r+2 b r^2+q\right)}{f^2 r^2},$$
$$R_{1323}= R_{1424}=-\frac{2 b r^2-q}{f^2 r}, \ R_{3434}=\frac{r \left(-2 a r+4 b r^2+2 q+F r\right)}{f^4}.$$
$$S_{11}=-\frac{8 b \left(-a r+2 b r^2+q\right)}{r^2}, \ S_{12}=-\frac{4 b}{r}, \ S_{33}= S_{44}=-\frac{2 a-8 b r-F}{f^2}$$
$$\mbox{and } \kappa = -\frac{2 (-2 a+12 b r+F)}{r^2},$$
where $F = f_3^2+f_4^2-f (f_{33}+f_{44})$.
Also the non-zero components of $C$ and $P$ are given by
$$C_{1212}=\frac{2 a r-6 q-F r}{3 r^3}, \ \ C_{1313}= C_{1414}= -\frac{\left(-a r+2 b r^2+q\right) (2 a r-6 q-F r)}{3 f^2 r^2}$$
$$C_{1323}= C_{1424}=-\frac{2 a r-6 q-F r}{6 f^2 r}, \ \ C_{3434}=-\frac{r (2 a r-6 q-F r)}{3 f^4},$$
$$-P_{1212}= P_{1221}=\frac{2 \left(3 q-2 b r^2\right)}{3 r^3}, \ P_{1313}= P_{1414}=-\frac{2 \left(2 b r^2-3 q\right) \left(-a r+2 b r^2+q\right)}{3 f^2 r^2},$$
$$P_{1323}= P_{1424}= P_{2313}= P_{2414}=\frac{3 q-2 b r^2}{3 f^2 r},$$
$$P_{1331}= P_{1441}=-\frac{2 \left(-a r+2 b r^2+q\right) \left(-2 a r+2 b r^2+3 q+F r\right)}{3 f^2 r^2},$$
$$P_{1332}= P_{1442}= P_{2331}= P_{2441}=-\frac{-2 a r+2 b r^2+3 q+F r}{3 f^2 r},$$
$$P_{3434}= -P_{3443}=\frac{2 r \left(-2 a r+2 b r^2+3 q+F r\right)}{3 f^4}.$$
The non-zero components of $\nabla R$ and $\nabla S$ are given by
$$R_{1212,2}=\frac{6 q}{r^4}, \ R_{1213,3}= R_{1214,4}= R_{1313,2}= R_{1414,2}=\frac{2 \left(2 b r^2-3 q\right) \left(-a r+2 b r^2+q\right)}{f^2 r^3},$$
$$R_{1223,3}= R_{1224,4}= R_{1323,2}= R_{1424,2}=\frac{2 b r^2-3 q}{f^2 r^2}, \ R_{3434,3}=\frac{F_3 r^2}{f^4},$$
$$R_{2334,4}= -R_{2434,3}= -\frac{1}{2}R_{3434,2}=\frac{-2 a r+2 b r^2+3 q+F r}{f^4}, \ R_{3434,4}=\frac{F_4 r^2}{f^4},$$
$$S_{11,2}=\frac{8 b \left(-a r+2 b r^2+q\right)}{r^3}, \ S_{33,3}= S_{44,3}=\frac{F_3}{f^2},\ S_{33,4}= S_{44,4}=\frac{F_4}{f^2},$$
$$S_{12,2}=\frac{4 b}{r^2}, \ 2S_{23,3}= 2S_{24,4}= S_{33,2}= S_{44,2}=-\frac{2 (-2 a+4 b r+F)}{f^2 r}.$$
The non-zero components of $R\cdot R$, $Q(g,R)$ and $Q(S,R)$ are given by
$$-R\cdot R_{121323}= -R\cdot R_{121424}= R\cdot R_{122313}= R\cdot R_{122414}=\frac{\left(2 b r^2-3 q\right) \left(2 b r^2-q\right)}{f^2 r^4},$$
$$R\cdot R_{133414}= -R\cdot R_{143413}=\frac{2 \left(2 b r^2-q\right) \left(-a r+2 b r^2+q\right) \left(-2 a r+2 b r^2+3 q+F r\right)}{f^4 r^3},$$
$$R\cdot R_{133424}= -R\cdot R_{143423}= R\cdot R_{233414}= -R\cdot R_{243413}=\frac{\left(2 b r^2-q\right) \left(-2 a r+2 b r^2+3 q+F r\right)}{f^4 r^2}.$$
$$Q(g,R)_{121323}= Q(g,R)_{121424}= -Q(g,R)_{122313}= -Q(g,R)_{122414}=\frac{2 b r^2-3 q}{f^2 r},$$
$$-Q(g,R)_{133414}= Q(g,R)_{143413}=\frac{2 \left(-a r+2 b r^2+q\right) \left(-2 a r+2 b r^2+3 q+F r\right)}{f^4},$$
$$-Q(g,R)_{133424}= Q(g,R)_{143423}= -Q(g,R)_{233414}= Q(g,R)_{243413}=\frac{r \left(-2 a r+2 b r^2+3 q+F r\right)}{f^4},$$
$$-Q(S,R)_{121323}= -Q(S,R)_{121424}= Q(S,R)_{122313}= Q(S,R)_{122414}=\frac{4 a q+8 b^2 r^3-20 b q r-2q F}{f^2 r^3},$$
$$-Q(S,R)_{133414}= Q(S,R)_{143413}=\frac{2 \left(-a r+2 b r^2+q\right) \left(4 a b r^2+2 a q-16 b q r-2 b F r^2-q F\right)}{f^4 r^2},$$
$$Q(S,R)_{133424}= Q(S,R)_{143423}= Q(S,R)_{233414}= Q(S,R)_{243413}$$
$$=-\frac{-4 a b r^2-2 a q+16 b q r+2 b F r^2+q F}{f^4 r}.$$
\section{\bf Conclusion}\label{con}
 From above we can conclude that the metric \eqref{RTII} fulfills the following curvature restricted geometric structures.
\begin{thm}
The Robinson-Trautman metric of Petrov type II given in \eqref{RTII} possesses the following curvature properties:
\begin{enumerate}[label=(\roman*)]
    \item It is Deszcz pseudosymmetric with $R\cdot R = \frac{q - 2 b r^2}{r^3} Q(g,R)$. Hence it is Ricci pseudosymmetric, conformally pseudosymmetric, projectively pseudosymmetric, concircularly pseudosymmetric and conharmonically pseudosymmetric of non-constant type with same associated scalar.
		\item  It satisfies
		$${r^3 \left(2 a r-6 q- F r\right)}\left[R\cdot R - Q(S, R)\right]= {2 \left(4 a q r+4 b^2 r^4-12 b q r^2-3 q^2-2 q F r\right)} Q(g,C).$$
		\item  It satisfies $C\cdot R = L Q(g, R)$, $L = -\frac{\left(2 a r-6 q - F r\right)}{6 r^3}$. Hence $C\cdot S = L Q(g, S)$, $C\cdot C = L Q(g, C)$ (i.e., pseudosymmetric Weyl conformal curvature tensor), $C\cdot W = L Q(g, W)$, $C\cdot P = L Q(g, P)$ and $C\cdot K = L Q(g, K)$.
		\item If $\left(4 q r (a-3 b r)+4 b^2 r^4-3 q^2-2 q F r\right)$ is nowhere zero, then it satisfies $C\cdot R - R\cdot C = L_1 Q(g,R) +L_2 Q(S,R)$ for
		$$L_1 = -\frac{(2 a r-6 q-F r) \left(2 \left(a q+4 b^2 r^3-6 b q r\right)-q F\right)}{3 r^2 \left(4 q r (a-3 b r)+4 b^2 r^4-3 q^2-2 q F r\right)}$$
		and 
		$$L_2 = \frac{\left(q-2 b r^2\right) (-2 a r+6 q+F r)}{8 q r (3 b r-a)-8 b^2 r^4+6 q^2+4 q F r}$$
		\item  It satisfies $C\cdot R - R\cdot C = L_1 Q(g,C) +L_2 Q(S,C)$ for
		$$L_1 = \frac{2 (-2 a+12 b r+F)}{3 r^2} \mbox{ and } L_2 = 1$$
		\item If $-2 a r+6 q+F r$ is everywhere nonzero, then its conformal curvature 2-forms are recurrent for the 1-form
		$$\Pi = \left\{0, \frac{F-2 a}{-2 a r+6 q+F r}, \frac{F_3 r}{-2 a r+6 q+F r}, \frac{F_4 r}{-2 a r+6 q+F r}\right\}.$$
		\item If $-2 a+4 b r+F$ is everywhere nonzero, then this metric is a Roter type with $R = N_1 (S\wedge S) + N_2 (g\wedge S) + N_3 (g\wedge g),$ where
		$$N_1 = \frac{r \left(2 a r-6 q-F r\right)}{2 \left(-2 a+4 b r+F\right)^2},$$
		$$N_2 = \frac{4 a b r^2+6 a q+8 b^2 r^3-F \left(2 b r^2+3 q\right)-36 b q r}{r \left(-2 a+4 b r+F\right)^2}$$
		$$N_3 = -\frac{1}{r^3}\left[\-\frac{4 b r \left(6 a q+8 b^2 r^3-24 b q r-3 q F\right)}{\left(-2 a+4 b r+F\right)^2}+q\right].$$
		\item This is 2-quasi Einstein as $S-\alpha g$ is of rank 2 for $\alpha = \frac{2 a-8 b r-F}{r^2}$ or $\alpha = -\frac{4 b}{r}$.
		\item This is generalized quasi Einstein in the sense of Chaki for
		$$\alpha = -\frac{-2 a+8 b r+F}{r^2}, \ \beta = \frac{4 b \left(-2 a+4 b r+F\right)}{\Phi_1^2 r},$$
		$$\gamma = -2 a+4 b r+F, \ \ \Pi =\left\{\frac{q-a r}{\Phi_1 r^3},\frac{1}{\Phi_1 r^2},0,0\right\}, \ \ \Phi= \left\{\Phi_1,0,0,0\right\},$$
		where $\Phi_1$ being arbitrary scalar, $||\Pi|| = -\frac{4 b}{\eta (1)^2 r^3}$ and $||\Phi|| = 0$.
		\item This is generalized quasi Einstein in the sense of De and Ghosh for
		$$\alpha = -\frac{4 b}{r}, \ \beta = \gamma = \frac{-2 a+4 b r+F}{r^2},$$
		$$\Pi =\left\{0,0,\frac{r}{f},0\right\}, \ \ \Phi= \left\{0,0,0,\frac{r}{f}\right\},$$
		where $||\Pi|| = ||\Phi|| = -1$ and $\Pi$ is orthogonal to $\Phi$.
		\item This is pseudo quasi Einstein for
		$$\alpha = -\frac{F}{2 r^2}, \ \beta = -\frac{4 (a-6 b r)}{r^2}, \ \gamma = 1, \ \Pi =\left\{0,0,0,\frac{r}{f}\right\}.$$
		\item A (0,2) tensor of the form
$$\left(
\begin{array}{cccc}
 t(1,1) & t(1,2) & 0 & 0 \\
 t(1,2) & t(2,2) & 0 & 0 \\
 0 & 0 & t(3,3) & t(3,4) \\
 0 & 0 & t(3,4) & t(4,4)
\end{array}
\right),$$
where $t(i,j)$ being arbitrary scalar, is $R$-compatible as well as $C$-compatible, $W$-compatible and $K$-compatible, where $t(i,j)$ being arbitrary scalar. So the Ricci tensor is $R$-compatible as well as $C$-compatible. Again the vector of the form $\{t_1,0,0,0\}$ or $\{0,t_2,0,0\}$ or $\{0,0,t_3,0\}$ or $\{0,0,0,t_4\}$, $t_i$ being arbitrary scalar, is $R$-compatible as well as $C$-compatible.
		\end{enumerate}
\end{thm}
%
\begin{rem}
From the value of the local components (presented in Section \ref{com}) of various tensors of the Robinson-Trautman metric \eqref{RTII}, we can easily conclude that the metric does not fulfill the following geometric structures:
\begin{enumerate}[label=(\roman*)]
\item not conformally symmetric and hence not locally symmetric, projectively symmetric, concircularly symmetric, conharmonically symmetric.
\item not conformally recurrent and hence not recurrent, projectively recurrent, concircularly recurrent, conharmonically recurrent.
\item not super generalized recurrent and hence not hyper generalized recurrent, weakly generalized recurrent.
\item not weakly symmetric for $R$, $C$, $P$, $W$ and $K$ and hence not Chaki pseudosymmetric for $R$, $C$, $P$, $W$ and $K$.
\item not weakly Ricci symmetric and hence not pseudo Ricci symmetric in the sense of Chaki and not Ricci symmetric.
\indent Scalar curvature is not constant and hence its Ricci tensor is not Codazzi type, not cyclic parallel.
\item $div R \ne 0$, $div C \ne 0$, $div P \ne 0$, $div W \ne 0$ and $div K \ne 0$.
\item not conformally semisymmetric and hence not semisymmetric, projectively semisymmetric, concircularly semisymmetric, conharmonically semisymmetric.
\item not Ricci generalized pseudosymmetric.
\item curvature 2-forms for $R$, $P$, $W$ and $K$ are not recurrent. The Ricci 1-forms are not recurrent.
\end{enumerate}
\end{rem}
\section{\bf Energy momentum tensor of Robinson-Trautman metric}\label{em-t}
From the Einstein's field equations, the energy momentum tensor, $T$ is given by 
$$T= \frac{c^4}{8\pi G}\left[S-\left(\frac{\kappa}{2}-\Lambda\right)g\right],$$
where $c=$ speed of light in vacuum, $G=$ gravitational constant and $\Lambda=$ cosmological constant. Thus the components of the energy momentum tensor is calculated by considering the components of Ricci tensor and scalar curvature. Now the non-zero components of $T$ (upto symmetry) are given by
$$T_{11}=\frac{c^4 \left(-a r+2 b r^2+q\right) \left(-2 a+8 b r+F+\Lambda  r^2\right)}{4 \pi  G r^3}$$
$$T_{12}=\frac{c^4 \left(-2 a+8 b r+F+\Lambda  r^2\right)}{8 \pi  G r^2}, \ T_{33}= T_{44}=-\frac{c^4 r (4 b+\Lambda  r)}{8 \pi  f^2 G}.$$
Then the non-zero components of covariant derivative of $T$ are given by
$$T_{11,2}=-\frac{c^4 \left(-a r+2 b r^2+q\right) (-2 a+4 b r+F)}{2 \pi  G r^4}, \ T_{11,3}=\frac{c^4 F_3 \left(-a r+2 b r^2+q\right)}{4 \pi  G r^3},$$
$$T_{11,4}=\frac{c^4 F_4 \left(-a r+2 b r^2+q\right)}{4 \pi  G r^3}, \ T_{12,2}=-\frac{c^4 (-2 a+4 b r+F)}{4 \pi  G r^3},$$
$$T_{12,3}=\frac{c^4 F_3}{8 \pi  G r^2}, \ T_{12,4}=\frac{c^4 F_4}{8 \pi  G r^2},$$
$$T_{23,3}= T_{24,4}=-\frac{c^4 (-2 a+4 b r+F)}{8 \pi  f^2 G r}, \ T_{33,2}= T_{44,2}=\frac{b c^4}{2 \pi  f^2 G}.$$
\indent We note that from the above it is easy to check that $div T = 0$ which implies that the metric \eqref{RTII} satisfies energy condition. However it is obvious that $T$ is not covariantly constant as well as it is not Codazzi type and cyclic parallel. In this section we would like to investigate the conditions under which $T$ is parallel, $T$ is Codazzi type and $T$ is cyclic parallel.\\
\indent Now from the values of the components of $\nabla T$, we get
$$T_{11,2} + T_{12,1} + T_{21,1}=-\frac{c^4 \left(-a r+2 b r^2+q\right) (-2 a+4 b r+F)}{2 \pi  G r^4},$$
$$T_{11,3} + T_{13,1} + T_{31,1}=\frac{c^4 F_3 \left(-a r+2 b r^2+q\right)}{4 \pi  G r^3},$$
$$T_{11,4} + T_{14,1} + T_{41,1}=\frac{c^4 F_4 \left(-a r+2 b r^2+q\right)}{4 \pi  G r^3},$$
$$T_{12,2} + T_{22,1} + T_{21,2}=-\frac{c^4 (-2 a+4 b r+F)}{2 \pi  G r^3},$$
$$T_{12,3} + T_{23,1} + T_{31,2}=\frac{c^4 F_3}{8 \pi  G r^2}, \ T_{12,4} + T_{24,1} + T_{41,2}=\frac{c^4 F_4}{8 \pi  G r^2},$$
$$T_{23,3} + T_{33,2} + T_{32,3}= T_{24,4} + T_{44,2} + T_{42,4}=-\frac{c^4 (-2 a+2 b r+F)}{4 \pi  f^2 G r}$$
and
$$T_{12,1} - T_{11,2} = \frac{c^4 \left(-a r+2 b r^2+q\right) (-2 a+4 b r+F)}{2 \pi  G r^4},$$
$$T_{11,3} - T_{13,1} = \frac{c^4 F_3 \left(-a r+2 b r^2+q\right)}{4 \pi  G r^3}, \ \ 
T_{11,4} - T_{14,1} = \frac{c^4 F_4 \left(-a r+2 b r^2+q\right)}{4 \pi  G r^3},$$
$$T_{12,3} - T_{13,2} = T_{21,3} - T_{23,1} = \frac{c^4 F_3}{8 \pi  G r^2}, \ \ T_{12,4} - T_{14,2} = T_{21,4} - T_{24,1} = \frac{c^4 F_4}{8 \pi  G r^2},$$
$$T_{22,1} - T_{21,2} = \frac{c^4 (-2 a+4 b r+F)}{4 \pi  G r^3}, \ \ T_{33,2} - T_{32,3} =  T_{44,2} - T_{42,4} = \frac{c^4 (-2 a+8 b r+F)}{8 \pi  f^2 G r}.$$
\indent Hence from above we see that if the energy-momentum tensor $T$ is cyclic parallel or Codazzi type or covariantly constant, then $b = 0$. Therefore $T$ is cyclic parallel or Codazzi type or covariantly constant if and only if $b=0$ and $F=2a$.
Now we can state the following:
\begin{thm}
For the Robinson-Trautman metric, the following conditions are equivalent:\\
(i) the  energy-momentum tensor is covariantly constant\\
(ii) the  energy-momentum tensor is Codazzi type\\
(iii) the  energy-momentum tensor is cyclic parallel\\
(iv) $b=0$ and $F=2a$.
\end{thm}
\begin{exm}
We consider $a=b=0$ and $f(x^3, x^4) = e^{x^3+x^4}$ in \eqref{RTII}. Then $F=0$ and therefore the energy-momentum tensor $T$ is covariantly constant and hence cyclic parallel and Codazzi type.
\end{exm}
\indent Now since the Robinson-Trautman metric \eqref{RTII} is Ricci pseudosymmetric and $T$ is linear combination of $S$ and $g$, and $R\cdot g = Q(g,g) = 0$, then $R\cdot T =\frac{q - 2 b r^2}{r^3} Q(g, T)$. Thus we can state the following:
\begin{thm}
The energy-momentum tensor of Robinson-Trautman metric is pseudosymmetric type, i.e., $R\cdot T = \frac{q - 2 b r^2}{r^3} Q(g, T)$.
\end{thm}
\section{\bf Som-Raychaudhuri metric and Robinson-Trautman metric}\label{sr-rt}
Som-Raychaudhuri \cite{SR68} spacetime is a G\"odel type metric  (see \cite{DHJKS14}, \cite{RT83}, \cite{SK16srs} and references therein) and in cylindrical coordinates $(t, r, z, \phi)$, its line element is given by
\be\label{srm}
ds^2 = dt^2 - (r^2-a^2r^4) d\phi^2 - dr^2 - dz^2 + 2 a r^2 d\phi dt.
\ee 
Recently, Shaikh and Kundu \cite{SK16srs} investigated the curvature restricted geometric structures admitting by the Som-Raychaudhuri spacetime and showed that such a spacetime is a 2-quasi-Einstein, generalized Roter type,  $Ein(3)$ manifold satisfying $R.R = Q(S,R)$, $C\cdot C = \frac{2a^2}{3} Q(g,C)$, and its Ricci tensor is cyclic parallel and Riemann compatible. In this section we study the comparison between the curvature properties of Som-Raychaudhuri metric and Robinson-Trautman metric.\\
\indent Physically, Som-Raychaudhuri \cite{SR68} spacetime is a stationary cylindrical symmetric solution of Einstein field equation corresponding to a charged dust distribution in rigid rotation whereas Robinson-Trautman metric admit a geodesic, hypersurface orthogonal, shear-free and expanding null congruences. Now from Section \ref{con} of this paper and Section 3 of \cite{SK16srs} we have the following comparison between Som-Raychaudhuri metric and Robinson-Trautman metric:\\
\textbf{A. Similarity:}\\
(i) Both are non-Einstein but 2-quasi-Einstein,\\
(ii) Both are pseudo quasi Einstein and generalized quasi-Einstein in the sense of Chaki as well as in the sense of De and Ghosh.\\
(iii) Both are of pseudosymmetric Weyl conformal curvature tensor.\\
(iv) Their Ricci tensors are Riemann compatible as well as conformal compatible, concircular compatible and conharmonic compatible.\\
\textbf{B. Dissimilarity:}\\
(i) Som-Raychaudhuri metric is Ricci generalized pseudosymmetric but Robinson-Trautman metric is Deszcz pseudosymmetric.\\
(ii) Som-Raychaudhuri metric is generalized Roter type but Robinson-Trautman metric is Roter type.\\
(iii) Ricci tensor of Som-Raychaudhuri metric is cyclic parallel whereas the scalar curvature of Robinson-Trautman metric is non-constant.\\
(iv) Robinson-Trautman metric is $Ein(2)$ but Som-Raychaudhuri metric is $Ein(3)$.


\end{document}